\documentclass[12pt]{article}
\usepackage{booktabs}
\usepackage[margin=2cm,a4paper]{geometry} 
\usepackage{amssymb, amsthm, amsmath}

\usepackage{pstricks,pst-node,pst-tree}

\usepackage{tikz}

\newcommand{\fra}[2]{\textstyle{\frac{#1}{#2}}}
\newcommand{\sgn}[1]{\text{sgn}(#1)}

\theoremstyle{plain}
\newtheorem{thm}{Theorem}[section]

\newtheorem{cor}{Corollary}
\newtheorem{res}{Result}

\newcommand{\lrr}[1]{{\langle #1 \rangle}_{\mathbb{R}}}
\newcommand{\id}{{\bf 1}}

\newcommand{\mat}[4]{\left(\begin{smallmatrix}#1 & #2 \\ #3 & #4 \end{smallmatrix}\right)}
\newcommand{\Mat}[4]{\left(\begin{matrix}#1 & #2 \\ #3 & #4 \end{matrix}\right)}
\newcommand{\lie}[2]{\left[#1 , #2\right]}
\newcommand{\beqn}{\begin{eqnarray}\begin{aligned}}
\newcommand{\eqn}{\end{aligned}\end{eqnarray}}

\begin{document}

\title{Lie geometry of $2\times 2$ Markov matrices}

\author{Jeremy G. Sumner\\
\small{School of Mathematics and Physics}\\
\small{University of Tasmania, Australia}\\
\small{ARC Research Fellow}\\
\small{\textit{email:} jsumner@utas.edu.au}}
\maketitle

\abstract{\noindent In recent work discussing model choice for continuous-time Markov chains, we have argued that it is important that the Markov matrices that define the model are closed under matrix multiplication \cite{sumner2012a,sumner2012b}.
The primary requirement is then that the associated set of rate matrices form a Lie algebra.
For the generic case, this connection to Lie theory seems to have first been made by \cite{johnson1985}, with applications for specific models given in \cite{bashford2004} and \cite{house2012}. 
Here we take a different perspective: given a model that forms a Lie algebra, we apply existing Lie theory to gain additional insight into the geometry of the associated Markov matrices.
In this short note, we present the simplest case possible of $2\times 2$ Markov matrices.
The main result is a novel decomposition of $2\times 2$ Markov matrices that parameterises the general Markov model as a perturbation away from the binary-symmetric model.
This alternative parameterisation provides a useful tool for visualising the binary-symmetric model as a submodel of the general Markov model.}

\vspace{2em}

\noindent \small{\textit{keywords:} Lie algebras, algebra, symmetry, Markov chains, phylogenetics}

\section{Results}

Consider the set of real $2\times 2$ Markov matrices
\beqn
\left\{\Mat{1-a}{b}{a}{1-b}:a,b\in \mathbb{R}\right\}\nonumber,
\eqn
and the subset of $2\times 2$ ``stochastic'' Markov matrices
\beqn
\left\{\Mat{1-a}{b}{a}{1-b}:0\leq a,b\leq 1\in \mathbb{R}\right\}\nonumber.
\eqn
In models of phylogenetic molecular evolution (see for example \cite{semple2003}), this set provides the transition matrices for what is known as the ``general Markov model'' on two states.
If we were to take the additional constraint $a\!=\!b$, the model would then be referred to as ``binary-symmetric''. 

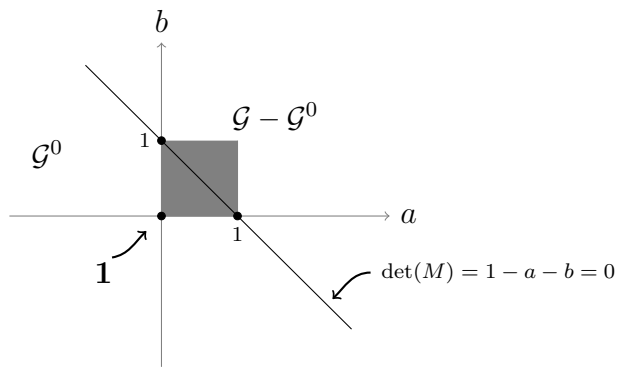
\begin{figure}[ht]
\centering

\begin{tikzpicture}{scale=4}
	\draw [fill,gray] (0,0) to (1,0) to (0,1);
	\draw [fill,gray] (1,1) to (1,0) to (0,1);
    \draw [thin, gray, ->] (0,-2) -- (0,2.3)      % draw y-axis line
        node [above, black] {$b$};              % add label for y-axis

    \draw [thin, gray, ->] (-2,0) -- (3,0)      % draw x-axis line
        node [right, black] {$a$};              % add label for x-axis
    \draw [] (2.5,-1.5) -- (-1,2);  % draw the graph

    \node [left] at (0,1) {\scriptsize {$1$}};                % label y-intercept
    \node [below] at (1,0) {\scriptsize {$1$}};               % label x-intercept
    \node [right] at (2.75,-.75) {\scriptsize{$\det(M)=1-a-b=0$}};
    \draw [->,thick] (2.75,-.75) to [out=180,in=45] (2.25,-1.1);
    \node [left] at (-.5,-.75) {$\id$};
    \node [above] at (-1.5,.5) {\small {$\mathcal{G}^0$}};
    \node [above] at (1.5,1) {\small {$\mathcal{G}-\mathcal{G}^0$}};
    \draw [<-,thick] (-.2,-.2) to [out=225,in=10] (-.65,-.55) ;
    \draw[fill] (0,1) circle [radius=0.05];
    \draw[fill] (1,0) circle [radius=0.05];
    \draw[fill] (0,0) circle [radius=0.05];
   
\end{tikzpicture}

\caption[Geometric rep]{The group $\mathcal{G}$ of invertible $2\times 2$ Markov matrices of the form $\mat{1-a}{b}{a}{1-b}$ understood geometrically as a manifold in $\mathbb{R}^2$. 
The gray area indicates the subset of ``stochastic'' Markov matrices.
The line $\det(M)=0$ indicates the boundary of the connected component to the identity, $\mathcal{G}^0$.}
\label{fig:mani}
\end{figure}

Associated with these sets is the matrix group 
\beqn
\mathcal{G}:=\left\{\Mat{1-a}{b}{a}{1-b}:a,b\in \mathbb{R},a+b\neq 1\right\}\nonumber.
\eqn
We can geometrically understand $\mathcal{G}$ by considering it as a manifold in $\mathbb{R}^2$.
This is illustrated in Figure~\ref{fig:mani}.

By considering smooth paths $A(t)\in \mathcal{G}$, we can define the tangent space of this matrix group at the identity:
\beqn
T_{\id}(\mathcal{G})=\left\{A'(0):A(t)\in \mathcal{G}\text{ and }A(0)=\id\right\}.\nonumber
\eqn
As $\mathcal{G}$ is a matrix group, it follows that $T_{\id}(\mathcal{G})$ forms a \emph{Lie algebra}.
This means that for all $X,Y\in T_{\id}(\mathcal{G})$ and $\lambda\in \mathbb{R}$, we have:
\begin{enumerate}
\item $X+\lambda Y\in T_{\id}(G)$, ie. $T_{\id}(G)$ is a vector space,
\item $\lie{X}{Y}:=XY-YX\in T_{\id}(\mathcal{G})$.
\end{enumerate}

Consider two smooth functions $a(t)$ and $b(t)$ satisfying $a(t)+b(t)\neq 1$ for all $t$, and $a(0)=b(0)=0$.
Define 
\[A(t)=\Mat{1-a(t)}{b(t)}{a(t)}{1-b(t)}.\]
Then, by construction, $A(t)$ is a smooth path in $\mathcal{G}$ and $A'(0)\in T_{\id}(\mathcal{G})$.
If we define $L_1:=\mat{-1}{0}{1}{0}$ and $L_2:=\mat{0}{1}{0}{-1}$, we have $A'(0)=a'(0)L_1+b'(0)L_2$,
so $T_{\id}(\mathcal{G})=\lrr{L_1,L_2}$ and $\{L_1,L_2\}$ is a basis for $T_{\id}(\mathcal{G})$.
It is straightforward to check that $\lie{L_1}{L_2}=L_1-L_2$, so we conclude that $T_{\id}(\mathcal{G})$ is indeed a Lie algebra.

Recall that a subgroup $H\leq G$ of a group is \emph{normal} if $ghg^{-1}\in H$ for all $h\in H$ and $g\in G$.
Also recall that the connected component to the identity $G^0$  is normal in G.
In our case, this becomes:
\begin{res}
$\mathcal{G}^0=\left\{M\in \mathcal{G}: \det(M)>0\right\}$.
\end{res}
\begin{proof}
Consider $M=\mat{1-a}{b}{a}{1-b}=e^{Qt}$ where $Q:=\mat{-\alpha}{\beta}{\alpha}{-\beta}$ is a rate matrix (as would occur in a continuous-time formulation of a Markov process).
Using the power series expansion of $e^{Qt}$, it is straightforward to show that, if $(a+b)t<1$,
\beqn
\alpha&=\frac{-\log(1-(a+b)t)}{1+b/a},\quad \beta&=\frac{-\log(1-(a+b)t)}{1+a/b},\nonumber
\eqn
provides a solution to $M=e^{Qt}$.
If we define the path $A(t):=e^{Qt}$, we have $A(0)=\id$ and $A(1)=M$.
Thus, $M\in \mathcal{G}^0$ for all $a+b<1$.
On the other hand, if $a+b\geq 1$, there can be no path $B(t)\in \mathcal{G}$ with $B(0)=\id$ and $B(1)=M$ because we would have $\det(B(\tau))=0$ for some $\tau$ in the interval $(0,1]$.
\end{proof}

\begin{cor}
$\mathcal{G}^0=\left\{e^{Q}:Q=\Mat{-\alpha}{\beta}{\alpha}{-\beta};\alpha,\beta\in \mathbb{R}\right\}$.
\end{cor}

Recall the homomorphism theorem for groups (see for example \cite{stillwell2008}), which ensures, for any group homomorphism $\rho:G\rightarrow G'$, that (i.) the image of $\rho$ is a subgroup of $G'$, (ii.) the kernal $K$ of $\rho$ is normal in $G$, and (iii.) $G/K\cong G'$. 
To understand the set difference $\mathcal{G}-\mathcal{G}^0$, we notice that $\mathcal{G}^0$ is the kernal of the homomorphism,
\beqn
\mathcal{G} &\rightarrow \{1,-1\}\cong \mathbb{Z}_2,\\
M &\mapsto \sgn{\det(M)}.\nonumber
\eqn
The kernal of this homomorphism is $\mathcal{G}^0$, thus $\mathcal{G}/\mathcal{G}^0 = \{\mathcal{G}^0,P\mathcal{G}^0\}\cong \mathbb{Z}_2$ for some $P\in \mathcal{G}-\mathcal{G}^0$.
For reasons of symmetry, we reflect the identity $\id$ in the line $\det(M)$ and set $P=\mat{0}{1}{1}{0}$, noting that $P^2=\id$.
As $\mathcal{G}/\mathcal{G}^0$ is a partition of $\mathcal{G}$, we see that 
$\mathcal{G}-\mathcal{G}^0=P\mathcal{G}^0$ and $\mathcal{G}=\mathcal{G}^0\cup P \mathcal{G}^0$.

Somewhat trivially:
\begin{res}
As manifolds, $\mathcal{G}^0\cong P \mathcal{G}^0$.
\end{res}
\begin{proof}
Clearly,
\beqn
P: \mathcal{G}^0 &\rightarrow P \mathcal{G}^0\\
M &\mapsto PM,\nonumber
\eqn
is a diffeomorphism because it maps continuous paths to continuous paths.
\end{proof}
In particular, this means that:
\begin{res}
$\mathcal{G}^0$ is connected $\Leftrightarrow$ $P \mathcal{G}^0$ is connected
\end{res}
\begin{proof}
No proof is required, but we give one regardless to illustrate.
Consider the path $A(t)=e^{Q_2t}e^{Q_1(1-t)}\in \mathcal{G}^0$ with $A(0)=M_1:=e^{Q_1}$ and $A(1)=M_2:=e^{Q_2}$.
Now, $B(t):=PA(t)$ is a path in $P \mathcal{G}^0$ with $B(0)=P M_1$ and $B(1)=P M_2$. 
As any two points in $P \mathcal{G}^0$ can be written in this way, we are done.
\end{proof}

Recall that the \emph{center} $Z(G)$ of a group $G$ is the set of all $g\in G$ such that $gh\!=\!hg$ for all $h\in G$.
In our case, suppose that $N=\mat{1-c}{d}{c}{1-d}\in Z(\mathcal{G})$. 
Setting $NM\!=\!MN$ implies:
\beqn
\Mat{1-c}{d}{c}{1-d}\Mat{1-a}{b}{a}{1-b}&=\Mat{\ast}{b(1-c-d)+d}{a(1-d-c)+c}{\ast}\\
&=\Mat{1-a}{b}{a}{1-b}\Mat{1-c}{d}{c}{1-d}\\
&=\Mat{\ast}{d(1-b-a)+b}{c(1-b-a)+a}{\ast},\nonumber
\eqn
which is true if and only if $-bc\!=\!-ad$ for all $a$ and $b$.
This can only happen if $c=d=0$, thus $Z(\mathcal{G})=\{\id\}$.
Now, consider the basic theorem (see for example \cite{stillwell2008}):
\begin{thm}
If a matrix group $G$ is path connected with discrete center, then any non-discrete normal subgroup $H$ will have tangent space $T_{\id}(H)\neq \{0\}$.
Further, $T_{\id}(H)$ is an ideal of $T_\id(G)$, ie. $[X,Y]\in T_{\id}(H)$ for all $X\in T_{\id}(H)$ and $Y\in T_{\id}(G)$.
Therefore, any such $H$ can be detected by checking for ideals of $T_\id(G)$.
\end{thm}

In our case, $\mathcal{G}^0$ satisfies the conditions of this theorem.
Suppose $\mathcal{I}$ is a proper ideal of $T_\id(\mathcal{G}^0)$.
Then $\mathcal{I}$ is one-dimensional, and $Y:=xL_1+yL_2\in \mathcal{I}$ satisfies:
\beqn
\lie{Y}{L_1}=y(L_2-L_1),\text{ and }\lie{Y}{L_2}=x(L_1-L_2)\in \mathcal{I},\nonumber
\eqn
which can only be true if $Y\propto (L_1-L_2)$.

\begin{res}
$\lrr{Y}=\lrr{L_1-L_2}$ is the only proper ideal of $T_\id(\mathcal{G}^0)$.
\end{res}

We take $Y\!=\!L_1-L_2$ and note that $Y^2=0$, so $e^{Ys}=e^{(L_1-L_2)s}=\id+Ys=\mat{1-s}{-s}{s}{1+s}:=h_s$.
If we define the matrix group $\mathcal{H}:=\left\{\mat{1-s}{-s}{s}{1+s},s\in \mathbb{R}\right\}$, it is easy to confirm that $\mathcal{H}$ is normal in $\mathcal{G}^0$ and has tangent space $T_\id(\mathcal{H})=\lrr{Y}$.

Let $\mathbb{R}^\times_{> 0}$ be the set of positive real numbers considered as a group under multiplication.
We have:
\begin{res}
$\mathcal{H}$ is the kernal of the homomorphism $\mathcal{G}^0\rightarrow \mathbb{R}^\times_{> 0}$ defined by $M\mapsto \det(M)$.
Thus $\mathcal{G}^0/\mathcal{H}\cong \mathbb{R}^\times_{> 0}$.
\end{res}
\begin{proof}
\beqn
\det(M)=1 \Leftrightarrow a+b=0 \Leftrightarrow M=\Mat{1-a}{-a}{a}{1+a}.\nonumber
\eqn
\end{proof}

Since $h_sh_t=h_{s+t}$, i.e. $\mathcal{H}$ forms a one-parameter subgroup of $\mathcal{G}^0$, we have $\mathcal{H}\cong \mathbb{R}^+$, where $\mathbb{R}^+=\mathbb{R}$ is considered as a group under addition.
Note that $\mathcal{G}^0/\mathcal{H}$ is a parameterised partition of $\mathcal{G}^0$, so we can write $\mathcal{G}^0/\mathcal{H}=\cup_{t\in \mathbb{R}}e^{Qt}\mathcal{H}$, where $Q\in T_\id(\mathcal{G}^0)-T_\id(\mathcal{H})$. 
We then see that any $M\in \mathcal{G}^0$ can be written as a product $e^{Qt}h_s$, where $\det(M)=\det(e^{Qt})$.
Again for reasons of symmetry, we take $Q=\fra{1}{2}\mat{-1}{1}{1}{-1}$, i.e. $Q$ is the generator of the binary-symmetric model, and we have $\det(e^{Qt})=e^{-t}:=\lambda$.

This brings us to our main result.
\begin{res}
Any $M\in \mathcal{G}^0$ can be expressed as
\beqn\label{eq:decomp}
M=\Mat{1-a}{b}{a}{1-b}=e^{Qt}h_s&=\fra{1}{2}\Mat{1+e^{-t}}{1-e^{-t}}{1-e^{-t}}{1+e^{-t}}\Mat{1-s}{-s}{s}{1+s}\\
&=\fra{1}{2}\Mat{1+\lambda}{1-\lambda}{1-\lambda}{1+\lambda}\Mat{1-s}{-s}{s}{1+s},
\eqn
where $\det(M)=\lambda=e^{-t}=1-a-b$, and $s=\fra{1}{2}(a-b)\det(M)^{-1}$.
\end{res}
For the binary-symmetric model implemented as a stationary Markov chain, the parameter $\lambda=e^{-t}$ is proportional to the expected number of transitions in chain in time $t$.
Therefore we can think of the parameter $s$ as providing a perturbation away from the binary-symmetric model.
Additionally, to ensure that $M$ is a stochastic Markov matrix, with $a,b\geq 0$, we require
$-\fra{1}{2}(e^{t}-1)\leq s \leq \fra{1}{2}(e^t-1)$.

%\begin{figure}
%\centering
%\begin{tikzpicture}{scale=3}
%	
%	     \draw[color=white,fill=lightgray,domain=0:2.7]  plot(\x,{(exp(\x)-1)}) -- (2.7,1) -- (2.7,0); % The top fill 
%    \draw[color=white,fill=lightgray,domain=0:2.7]  plot(\x,{(exp(\x)-1)}) -- (2.7,-1) -- (2.7,0); % The bottom fill
%    \draw[domain=0:3.2]  plot(\x,{(exp(\x)-1)}) ; % The top graph
%    \draw[domain=0:3.2]  plot(\x,{(exp(\x)-1)}) ; % The bottom graph
%    
%    \draw [thin, gray, ->] (0,-2) -- (0,2.3)      % draw t-axis line
%        node [above, black] {$s$};              % add label for t-axis
%
%    \draw [thin, gray, ->] (-2,0) -- (4,0)      % draw s-axis line
%        node [right, black] {$t$};              % add label for s-axis
%    %\draw[fill] (0,0) circle [radius=0.05];
%    \draw[fill] (0,-1) circle [radius=0.05];
%    \draw[fill] (0,1) circle [radius=0.05];
%
%    \draw [dashed] (0,1) -- (3.5,1);
%     \draw [dashed] (0,-1) -- (3.5,-1);
%     
%     
%     \node [left] at (0,1) {\scriptsize{$\fra{1}{2}$}};
%     \node [left] at (0,-1) {\scriptsize{$-\fra{1}{2}$}};
%
%\end{tikzpicture}
%\caption{The region of stochastic Markov matrices represented under the parameterisation $s$ and $t$ provided in (\ref{eq:decomp}).}
%\label{fig:stochreg}
%\end{figure}

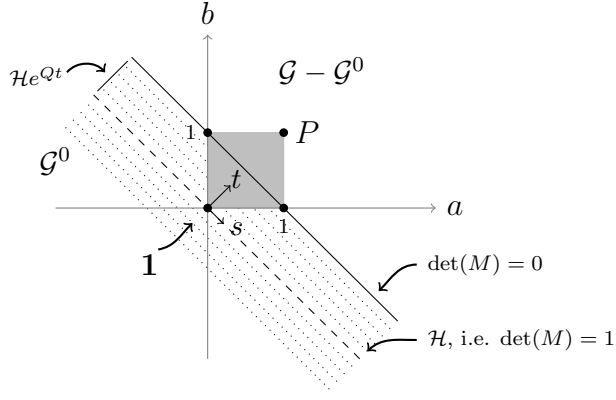
\begin{figure}
\centering
\begin{tikzpicture}{scale=3}

	\draw [dotted] (2.1,-1.9) -- (-1.4,1.6); % strata
    \draw [dotted] (2.2,-1.8) -- (-1.3,1.7); 
    \draw [dotted] (2.3,-1.7) -- (-1.2,1.8); 
    \draw [dotted] (2.4,-1.6) -- (-1.1,1.9); 
    
    \draw [dotted] (1.9,-2.1) -- (-1.6,1.4); % more strata
    \draw [dotted] (1.8,-2.2) -- (-1.7,1.3); 
    \draw [dotted] (1.7,-2.3) -- (-1.8,1.2); 
    \draw [dotted] (1.6,-2.4) -- (-1.9,1.1);

	\draw [fill=lightgray,lightgray] (0,0) to (1,0) to (1,1) to (0,1);
    \draw [thin, gray, ->] (0,-2) -- (0,2.3)      % draw y-axis line
        node [above, black] {$b$};              % add label for y-axis

    \draw [thin, gray, ->] (-2,0) -- (3,0)      % draw x-axis line
        node [right, black] {$a$};              % add label for x-axis
    \draw [] (2.5,-1.5) -- (-1,2);  % det(M)=0
    \draw [dashed] (2,-2) -- (-1.5,1.5);  % det(M=1)

    \node [left] at (0,1) {\scriptsize {$1$}};                	% label y-intercept
    \node [below] at (1,0) {\scriptsize {$1$}};               	% label x-intercept
    
    \node [right] at (2.75,-.75) {\scriptsize{$\det(M)=0$}};  	% det(M)=0
    \draw [->,thick] (2.75,-.75) to [out=180,in=45] (2.25,-1.1); 		% arrow from det(M)=0
    \node [right] at (2.75,-1.75) {\scriptsize{$\mathcal{H}$, i.e. $\det(M)=1$}};	% det(M)=1
    \draw [->,thick] (2.75,-1.75) to [out=180,in=45] (2.1,-1.9); 					% arrow from det(M)=1
    
	\node [right] at (-2.75,1.75) {\scriptsize{$\mathcal{H}e^{Qt}$}}; 
	\draw [-] (-1.45,1.55) to  (-1.05,1.95);
    \draw [<-,thick] (-1.35,1.8) to [out=130,in=50] (-1.85,1.8); 					%

    \node [left] at (-.5,-.75) {$\id$};
    \node [above] at (-2,.3) {\small {$\mathcal{G}^0$}};
    \node [above] at (1.5,1.5) {\small {$\mathcal{G}-\mathcal{G}^0$}};
    \draw[fill] (1,1) circle [radius=0.05];
    \node [right] at (1,1) {$P$};
    
    \draw [<-,thick] (-.2,-.2) to [out=225,in=10] (-.65,-.55) ;
    \draw[fill] (0,1) circle [radius=0.05];
    \draw[fill] (1,0) circle [radius=0.05];
    \draw[fill] (0,0) circle [radius=0.05];
    
	\draw[->] (0,0) -- (.3,.3); % arrow for t
    \node[] at (.38,.38) {\small{$t$}};

    \draw[->] (0,0) -- (.21,-.21); % arrow for s
    \node[] at (.38,-.28) {\small{$s$}};

\end{tikzpicture}
\caption{Lie geometry of $2\times 2$ Markov matrices}
\label{fig:liegeo}
\end{figure}

The decomposition (\ref{eq:decomp}) is the main result of this note and is presented geometrically in Figure~\ref{fig:liegeo}.
It is remarkable that such a simple application of elementary Lie theory has led directly to this decomposition, and it seems plausible that this decomposition may be useful in practice for (i.) computational efficiency, and/or (ii.) the simple interpretation of the parameters $t$, $\lambda$ and $s$.
It will be interesting to explore whether a similar analysis leads to alternative parameterisation for other popular phylogenetic models that form Lie algebras, but we leave this for future work.

\bibliographystyle{plain}
\bibliography{masterAB}

\begin{thebibliography}{1}

\bibitem{bashford2004}
J.~D. Bashford, P.~D. Jarvis, J.~G. Sumner, and M.~A. Steel.
\newblock {$U(1)\times U(1)\times U(1)$} symmetry of the {K}imura 3{S}{T} model
  and phylogenetic branching processes.
\newblock {\em J. Phys. A Math. Gen.}, 37:L1--L9, 2004.

\bibitem{house2012}
T.~House.
\newblock Lie algebra solution of population models based on time-inhomogeneous
  markov chains.
\newblock {\em J. Appl. Prob.}, 49:472--481, 2012.

\bibitem{johnson1985}
J.~E. Johnson.
\newblock {M}arkov-type {Lie} groups in {$GL(n,{R})$}.
\newblock {\em J. Math. Phys.}, 26:252--257, 1985.

\bibitem{semple2003}
C.~Semple and M.~Steel.
\newblock {\em Phylogenetics}.
\newblock Oxford Press, 2003.

\bibitem{stillwell2008}
John Stillwell.
\newblock {\em Naive Lie Theory}.
\newblock Undergraduate Texts in Mathematics. Springer, New York, 2008.

\bibitem{sumner2012a}
J.~G. Sumner, J.~Fern\'andez-S\'anchez, and P.~D. Jarvis.
\newblock Lie markov models.
\newblock {\em J. Theor. Biol.}, 298:16--31, 2012.

\bibitem{sumner2012b}
J.~G. Sumner, P.~D. Jarvis, J.~Fern\'andez-S\'anchez, B.~T. Kaine, Michael~D.
  Woodhams, and B.~R Holland.
\newblock Is the general time-reversible model bad for molecular phylogenetics?
\newblock {\em Syst. Biol.}, 61:1069--1074, 2012.

\end{thebibliography}

\end{document}